\documentclass[preprint,12pt]{elsarticle}

\usepackage{amssymb}
\usepackage{Tr_95}
\newtheorem{thm}{Theorem}

\newdefinition{dfn}{Definition}
\newdefinition{exm}{Example}
\newproof{pf}{Proof}
\newtheorem{cor}{Corollary}
\newtheorem{lm}{Lemma}

\journal{...}

\begin{document}

\begin{frontmatter}



%\title{GROUP $SU$-ACTION AND ITS APPLICATIONS TO GROUP THEORY}
\title{ On some new difference sequence spaces of fractional order}

%% \address[label1]

\author[label1]{Serkan Demiriz\corref{cor1}}
\ead{serkan.demiriz@gop.edu.tr}\cortext[cor1]{Corresponding Author
(Tel: +90 356 252 16 16, Fax: +90 356 252 15 85)}
\author[label3]{Osman Duyar}
\ead{osman-duyar@hotmail.com}

\address[label1]{Department of Mathematics, Faculty of Arts and Science, Gaziosmanpa\c{s}a University,\\
 60250 Tokat, Turkey }
\address[label3]{Anatolian High School, 60200 Tokat, Turkey }
%----------------------------------------------------------------------

\begin{abstract}
Let $\Delta^{(\alpha)}$ denote the fractional difference operator.
In this paper, we define new difference sequence spaces
$c_0(\Gamma,\Delta^{(\alpha)},u)$ and
$c(\Gamma,\Delta^{(\alpha)},u)$. Also,  the $\beta-$ dual of the
spaces  $c_0(\Gamma,\Delta^{(\alpha)},u)$ and
$c(\Gamma,\Delta^{(\alpha)},u)$  are determined and calculated their
Schauder basis. Furthermore, the classes
$(\mu(\Gamma,\Delta^{(\alpha)},u):\lambda)$ where $\mu\in
\{c_{0},c\}$ and $\lambda \in \{c_{0},c,\ell_{\infty},\ell_1\}$ .
\end{abstract}


\end{frontmatter}

% \linenumbers

%------------------------INTRODUCTION---------------------------
%\section{Introduction}
%---------------------------------------------------------------
\noindent
%------------------------PRELIMINARIES---------------------------
\section{Preliminaries,background and notation}
%---------------------------------------------------------------
By a \textit{sequence space}, we mean any vector subspace of
$\omega$, the space of all real or complex valued sequences
$x=(x_k)$. The well-known sequence spaces that we shall use
throughout this paper are as following:

$\ell_{\infty}$: the space of all bounded sequences,

$c$: the space of all convergent sequences,

$c_{0}$: the space of all null sequences,

$cs$: the space of all sequences which form convergent series,

$\ell_1$: the space of all sequences which form absolutely
convergent series,

$\ell_p$: the space of all sequences which form $p$-absolutely convergent series,\\
where $1< p<\infty$.

%-----------------------------------------------------------------

Let $X,Y$ be two sequence spaces and $A=(a_{nk})$ be an infinite
matrix of real or complex numbers $a_{nk}$, where $n,k\in
\mathbb{N}$. Then, we say that $A$ defines a matrix mapping from $X$
into $Y$, and we denote it by writing $A:X\rightarrow Y$, if for
every sequence $x=(x_{k})\in \lambda$ the sequence
$Ax=\{(Ax)_{n}\}$, the $A$-transform of $x$, is in $Y$; where
\begin{equation}\label{1}
(Ax)_{n}=\sum_{k} a_{nk}x_{k}, \quad (n\in \mathbb{N}).
\end{equation}
For simplicity in notation, here and in what follows, the summation
without limits runs from $0$ to $\infty$. By $(X:Y)$, we denote the
class of all matrices $A$ such that $A:X\rightarrow Y$. Thus, $A\in
(X:Y)$ if and only if the series on the right side of (\ref{1})
converges for each $n\in \mathbb{N}$ and every $x\in X$, and we have
$Ax=\{(Ax)_{n}\}_{n\in \mathbb{N}}\in Y$ for all $x\in X$. A
sequence $x$ is said to be $A$-summable to $\alpha$ if $Ax$
converges to $\alpha$ which is called as the $A$-limit of $x$.
%-------------------------------------------------------------------

If a normed sequence space $X$ contains a sequence $(b_{n})$ with
the property that for every $x\in X$ there is a unique sequence of
scalars $(\alpha_{n})$ such that
$$
\lim_{n\rightarrow \infty}
\|x-(\alpha_{0}b_{0}+\alpha_{1}b_{1}+...+\alpha_{n}b_{n})\|=0,
$$
then $(b_{n})$ is called a \emph{Schauder basis} (or briefly
\emph{basis}) for $X$. The series $\sum \alpha_{k}b_{k}$ which has
the sum $x$ is then called the expansion of $x$ with respect to
$(b_{n})$, and written as $x=\sum \alpha_{k}b_{k}$.

%-------------------------------------------------------------------

A matrix $A=(a_{nk})$ is called a \emph{triangle} if $a_{nk}=0$ for
$k>n$ and $a_{nn}\neq 0$ for all $n\in \mathbb{N}$. It is trivial
that $A(Bx)=(AB)x$ holds for triangle matrices $A,B$ and a sequence
$x$. Further, a triangle matrix $U$ uniquely has an inverse
$U^{-1}=V$ that is also a triangle matrix. Then, $x=U(Vx)=V(Ux)$
holds for all $x\in \omega$. We write additionally $\mathcal{U}$ for
the set of all sequences $u$ such that $u_k\neq0$ for all  $k \in
\mathbb{N}$.

For a sequence space $X$, the \emph{matrix domain} $X_{A}$ of an
infinite matrix $A$ is defined by
\begin{equation}\label{2}
X_{A}=\left\{x=(x_{k})\in \omega: Ax\in \lambda\right\},
\end{equation}
which is a sequence space.  The approach constructing a new sequence
space by means of the matrix domain of a particular limitation
method has been recently employed by Wang \cite{w}, Ng and Lee
\cite{nglee}, Ayd{\i}n and Ba\c{s}ar \cite{cafb} and Altay and
Ba\c{s}ar \cite{bafb5}.
%They introduced the sequence
%spaces $(\ell_{\infty})_{N_{q}}$ and $c_{N_{q}}$ in \cite{w},
%$(\ell_{p})_{C_{1}}=X_{p}$ in \cite{nglee},
%$(c_{0})_{A^{r}}=a_{0}^{r}$ and $c_{A^{r}}=a_{c}^{r}$ in
%\cite{cafb1}, $(c_{0})_{E^{r}}=e_{0}^{r}$ and $c_{E^{r}}=e_{c}^{r}$
%in \cite{bafb1}, $(\ell_{p})_{E^{r}}=e_{p}^{r}$ and
%$(\ell_{\infty})_{E^{r}}=e_{\infty}^{r}$ in \cite{bafbmm}; where
%$1\leq p<\infty$ and $N_{q}, C_{1}$ and $ E^{r}$ denote the
%N\"{o}rlund means, Ces\`{a}ro means of order $1$, Euler means of
%order $r$, respectively, where $0<r<1$.

%--------------------------------------------------------------------
The gamma function may be regarded as a generalization of $n!$
($n-$factorial), where $n$ is any positive integer.  The gamma
function  $\Gamma$ is defined for all $p$ real numbers except the
negative integers and zero. It can be expressed as an improper
integral as follows:
\begin{equation}\label{gamma}
\Gamma(p)=\int_0^\infty e^{-t}t^{p-1}dt.
\end{equation}

From the equality  (\ref{gamma}) we deduce following properties:\\

(i) If $n\in\mathbb{N}$ then we have $\Gamma(n+1)=n!$

(ii) If $n\in\mathbb{R}-\{0,-1,-2,-3...\}$ then we have
$\Gamma(n+1)=n\Gamma(n)$.\\

For a proper fraction $\alpha$, Baliarsingh and Dutta have defined a
fractional difference operators $\Delta^{\alpha}:w\rightarrow w$,
$\Delta^{(\alpha)}:w\rightarrow w$ and their inverse in \cite{bail3}
as follows:
\begin{equation}\label{op}
\Delta^{\alpha}(x_{k})=\sum_{i=0}^\infty(-1)^i\frac{\Gamma(\alpha+1)}{i!\Gamma(\alpha+1-i)}x_{k+i}
\end{equation}

\begin{equation}\label{op2}
\Delta^{(\alpha)}(x_{k})=\sum_{i=0}^\infty(-1)^i\frac{\Gamma(\alpha+1)}{i!\Gamma(\alpha+1-i)}x_{k-i}
\end{equation}

\begin{equation}\label{op3}
\Delta^{-\alpha}(x_{k})=\sum_{i=0}^\infty(-1)^i\frac{\Gamma(1-\alpha)}{i!\Gamma(1-\alpha-i)}x_{k+i}
\end{equation}
and
\begin{equation}\label{op4}
\Delta^{(-\alpha)}(x_{k})=\sum_{i=0}^\infty(-1)^i\frac{\Gamma(1-\alpha)}{i!\Gamma(1-\alpha-i)}x_{k-i}
\end{equation}
where we assume throughout the series defined in
(\ref{op})-(\ref{op4}) are convergent. In particular, for
$\alpha=\frac{1}{2}$,

\begin{eqnarray}\label{space1}
\nonumber
\Delta^{1/2}x_k&=&x_{k}-\frac{1}{2}x_{k+1}-\frac{1}{8}x_{k+2}-\frac{1}{16}x_{k+3}-\frac{5}{128}x_{k+4}
-\frac{7}{256}x_{k+5}-...\\ \nonumber
\Delta^{-1/2}x_k&=&x_{k}+\frac{1}{2}x_{k+1}+\frac{3}{8}x_{k+2}+\frac{5}{16}x_{k+3}+\frac{35}{128}x_{k+4}
+\frac{63}{256}x_{k+5}+...\\ \nonumber
\Delta^{(1/2)}x_k&=&x_{k}-\frac{1}{2}x_{k-1}-\frac{1}{8}x_{k-2}-\frac{1}{16}x_{k-3}-\frac{5}{128}x_{k-4}
-\frac{7}{256}x_{k-5}-...\\ \nonumber
\Delta^{(-1/2)}x_k&=&x_{k}+\frac{1}{2}x_{k-1}+\frac{3}{8}x_{k-2}+\frac{5}{16}x_{k-3}+\frac{35}{128}x_{k-4}
+\frac{63}{256}x_{k-5}+....
\end{eqnarray}

Baliarsingh have been defined the spaces
$X(\Gamma,\Delta^{\alpha},u)$ for $X\in\{\ell_\infty,c_0,c\}$ by
introducing  the fractional difference operator $\Delta^{\alpha}$
and a positive fraction $\alpha$ in \cite{bail1}. In this article,
Baliarsingh have been studied some topological properties of the
spaces $X(\Gamma,\Delta^{\alpha},u)$ and established their
$\alpha-$, $\beta-$ and $\gamma-$duals.

%-------------------------------------------------------------------

Following \cite{bail1}, we introduce the sequence spaces
$c_0(\Gamma,\Delta^{(\alpha)},u)$ and
$c(\Gamma,\Delta^{(\alpha)},u)$ and obtain some results related to
these sequence spaces. Furthermore, we compute the $\beta-$dual of
the spaces  $c_0(\Gamma,\Delta^{(\alpha)},u)$ and
$c(\Gamma,\Delta^{(\alpha)},u)$ . Finally, we characterize some
matrix transformations on new sequence spaces.
%-------------------------------------------------------------------

\section{The sequence spaces $c_0(\Gamma,\Delta^{(\alpha)},u)$ and $c(\Gamma,\Delta^{(\alpha)},u)$ }

In this section, we define the sequence spaces
$c_0(\Gamma,\Delta^{(\alpha)},u)$ and
$c(\Gamma,\Delta^{(\alpha)},u)$ and examine the some topological
properties of this sequence spaces.
%-------------------------------------------------------------------

The notion of difference sequence spaces was introduced by Kýzmaz
\cite{kzmaz}. It was generalized by Et and Çolak \cite{metolak} as
follows:

Let $m$ be a non-negative integer. Then
$$
\Delta^{m}(X)=\{x=(x_k): \Delta^{m}x_k\in X\}
$$
where $\Delta^{0}x=(x_k)$,
$\Delta^{m}x=(\Delta^{m-1}x_k-\Delta^{m-1}x_{k+1})$ for all
$k\in\mathbb{N}$ and

\begin{equation}\label{delta}
\Delta^{m}x_k=\sum_{i=0}^{m}(-1)^i\left(m \atop i\right)x_{k+i}.
\end{equation}

Furthermore, Malkowsky E., et al.\cite{emmmss} have been introduced
the spaces
\begin{equation}\label{delta}
\Delta_u^{(m)}X=\{x\in \omega: \Delta_u^{(m)}x\in X\}
\end{equation}
where $\Delta_u^{(m)}x=u\Delta^{(m)}x$ for all $x\in\omega$ . In
this study, the operator $\Delta_u^{(m)}:\omega\rightarrow\omega$
was defined as follows:
\begin{equation}\label{delta2}
\Delta^{(m)}x_k=\sum_{i=0}^{m}(-1)^i\left(m \atop i\right)x_{k-i}.
\end{equation}

Let $\alpha$ be a proper fraction and $u\in\mathcal{U}$. We define
the sequence spaces $c_0(\Gamma,\Delta^{(\alpha)},u)$ and
$c(\Gamma,\Delta^{(\alpha)},u)$ as follows:

\begin{equation}\label{space1}
c_0(\Gamma,\Delta^{(\alpha)},u)=\{x\in\omega:
(\sum_{j=0}^ku_j\Delta^{(\alpha)}x_j)\in c_0\}
\end{equation}
and
\begin{equation}\label{space2}
c(\Gamma,\Delta^{(\alpha)},u)=\{x\in\omega:
(\sum_{j=0}^ku_j\Delta^{(\alpha)}x_j)\in c\}.
\end{equation}

Now, we define the triangle matrix
$\Delta_u^{(\alpha)}(\Gamma)=(\tau_{nk})$,
\begin{equation}\label{matrix}
\tau_{nk}=\left\{\begin{array}{ll}
  \displaystyle
  \sum_{i=0}^{n-k}(-1)^i\frac{\Gamma(\alpha+1)}{i!\Gamma(\alpha+1-i)}u_{i+k},
   & (0\leq k\leq n)\\ \displaystyle 0, & k>n
\end{array}\right.
\end{equation}
for all $k,n\in \mathbb{N}$. Further, for any sequence $x=(x_k)$ we
define the sequence $y=(y_k)$ which will be used, as the
$\Delta_u^{(\alpha)}(\Gamma)-$transform of $x$, that is
\begin{eqnarray}\label{space1}
\nonumber
y_k&=&\sum_{j=0}^ku_j\Delta^{(\alpha)}x_j~=~\sum_{j=0}^ku_j(x_{j}-\alpha
x_{j-1}+\frac{\alpha(\alpha-1)}{2!}x_{j-3}+...)\\
&=&\sum_{j=0}^{k}\bigg(\sum_{i=0}^{k-j}(-1)^i\frac{\Gamma(\alpha+1)}{i!\Gamma(\alpha+1-i)}u_{i+j}\bigg)x_{j}
\end{eqnarray}
for all $k\in\mathbb{N}$.

%Furthermore
%\begin{equation}\label{ters}
%x_k=\sum_{i=0}^{\infty}(-1)^i\frac{\Gamma(1-\alpha)}{i!\Gamma(1-\alpha-i)}\frac{y_{k-i}-y_{k-i-1}}{u_{k-i}}
%\end{equation}
%or if we define the $\left(\alpha \atop i\right)$ as follows
%$$
%\left(\alpha \atop
%i\right)=\frac{\alpha(\alpha-1)(\alpha-2)...(\alpha-i+1)}{i!}~~and~~\left(\alpha
%\atop 0\right)=1
%$$
%for $\alpha\in\mathbb{R}$ and $i\in\mathbb{N}$. Then we can obtain
%\begin{equation}\label{ters1}
%x_k=\sum_{i=0}^{k}(-1)^i\left(-\alpha \atop
%i\right)\frac{y_{k-i}-y_{k-i-1}}{u_{k-i}}=\sum_{i=0}^{k}(-1)^i\left({-\alpha}
%\atop {k-i}\right)\frac{y_{i}-y_{i-1}}{u_{i}}.
%\end{equation}
%.....\\

 It is natural that the spaces
$c_0(\Gamma,\Delta^{(\alpha)},u)$ and
$c(\Gamma,\Delta^{(\alpha)},u)$ may also be defined with the
notation of (\ref{2}) that
\begin{equation}\label{fac1}
c_0(\Gamma,\Delta^{(\alpha)},u)=(c_0)_{\Delta_u^{(\alpha)}(\Gamma)}~~
\textrm{and}~~~~c(\Gamma,\Delta^{(\alpha)},u)=c_{\Delta_u^{(\alpha)}(\Gamma)}.
\end{equation}
Before the main result let us give some lemmas, which we use
frequently throughout this study, with respect to
$\Delta^{(\alpha)}$ operator.

\begin{lm}\label{lem1}\cite[Theorem 2.2]{bail3}
$$\Delta^{(\alpha)}o\Delta^{(\beta)}=\Delta^{(\beta)}o\Delta^{(\alpha)}=\Delta^{(\alpha+\beta)}.$$
\end{lm}

\begin{lm}\label{lem2}\cite[Theorem 2.3]{bail3}

$$\Delta^{(\alpha)}o\Delta^{(-\alpha)}=\Delta^{(-\alpha)}o\Delta^{(\alpha)}=Id$$
where $Id$  the identity operator on $\omega$.
\end{lm}

\begin{thm}
The sequence spaces $c_0(\Gamma,\Delta^{(\alpha)},u)$ and
$c(\Gamma,\Delta^{(\alpha)},u)$  are $BK-$spaces  with the norm
\begin{equation}\label{norm}
\|x\|_{c_0(\Gamma,\Delta^{(\alpha)},u)}=\|x\|_{c(\Gamma,\Delta^{(\alpha)},u)}
=\sup_k\bigg|\sum_{j=0}^ku_j\Delta^{(\alpha)}x_j\bigg|.
\end{equation}

\end{thm}
\begin{pf}
Since (\ref{fac1}) holds and $c_{0},c$ are $BK-$ spaces with respect
to their natural norms (see \cite[pp. 16-17]{fb}) and the matrix
$\Delta_u^{(\alpha)}(\Gamma)=(\tau_{nk})$ is a triangle, Theorem
4.3.12 Wilansky \cite[pp. 63]{aw} gives the fact that
$c_0(\Gamma,\Delta^{(\alpha)},u)$ and
$c(\Gamma,\Delta^{(\alpha)},u)$ are $BK-$ spaces with the given
norms. This completes the proof .
\end{pf}

Now, we may give the following theorem concerning the isomorphism
between the spaces $c_0(\Gamma,\Delta^{(\alpha)},u)$,
$c(\Gamma,\Delta^{(\alpha)},u)$ and $c_0,c$, respectively:

\begin{thm}
The sequence spaces $c_0(\Gamma,\Delta^{(\alpha)},u)$ and
$c(\Gamma,\Delta^{(\alpha)},u)$  are linearly isomorphic to the
spaces $c_0$ and $c$, respectively, i.e,
$c_0(\Gamma,\Delta^{(\alpha)},u)\cong c_0$ and
$c(\Gamma,\Delta^{(\alpha)},u)\cong c$.
\end{thm}
\begin{pf}
We prove the theorem for the space $c(\Gamma,\Delta^{(\alpha)},u)$.
To prove this, we should show the existence of a linear bijection
between the spaces $c(\Gamma,\Delta^{(\alpha)},u)$ and $c$. Consider
the transformation $T$ defined, with the notation of (\ref{space1}),
from $c(\Gamma,\Delta^{(\alpha)},u)$ to $c$ by $x\mapsto
y=Tx=\Delta_u^{(\alpha)}(\Gamma)x$. The linearity of $T$ is clear.
Further, it is trivial that $x=\theta$ whenever $Tx=\theta$ and
hence $T$ is injective. Let be $y=(y_k)\in c$ and we define a
sequence $x=(x_k)\in c(\Gamma,\Delta^{(\alpha)},u)$ by

\begin{equation}\label{ters}
x_k=\sum_{i=0}^{\infty}(-1)^i\frac{\Gamma(1-\alpha)}{i!\Gamma(1-\alpha-i)}\frac{y_{k-i}-y_{k-i-1}}{u_{k-i}}.
\end{equation}
Then by Lemma \ref{lem2}, we deduce that
\begin{eqnarray}\label{space5}
\nonumber
\sum_{j=0}^ku_j\Delta^{(\alpha)}x_j&=&\sum_{j=0}^ku_j\Delta^{(\alpha)}\bigg(\sum_{i=0}^{\infty}(-1)^i
\frac{\Gamma(1-\alpha)}{i!\Gamma(1-\alpha-i)}\frac{y_{j-i}-y_{j-i-1}}{u_{j-i}}\bigg)\\
\nonumber
&=&\sum_{j=0}^ku_j\Delta^{(\alpha)}\bigg(\Delta^{(-\alpha)}\bigg(\frac{y_{j}-y_{j-1}}{u_{j}}\bigg)\bigg)\\
&=&\sum_{j=0}^k(y_{j}-y_{j-1})~=~y_k
\end{eqnarray}
Hence, $x\in c(\Gamma,\Delta^{(\alpha)},u)$ so $T$ is surjective.
Furthermore one can easily show that $T$ is norm preserving. This
complete the proof.
\end{pf}

%-------------------------------------------------------------------

%-------------------------------------------------------------------

%%-------------------------------------------------------------------
%\begin{lm}\label{l1}
%The Banach space $f$ has no Schauder basis.
%\end{lm}
%
%%-------------------------------------------------------------------
%
%Since, it is known that the matrix domain $\lambda_{A}$ of a normed
%sequence space $\lambda$ has a basis if and only if $\lambda$ has a
%basis whenever $A=(a_{nk})$ is a triangle (cf. \cite{ajem}) and the
%space $f$ has no Schauder basis by Lemma \ref{l1}, we have:
%\begin{cor}
%The space $g^{f}$ has no Schauder basis.
%\end{cor}
%
%%-------------------------------------------------------------------
\section{The $\beta$-Dual of The Spaces $c_0(\Gamma,\Delta^{(\alpha)},u)$ and
$c(\Gamma,\Delta^{(\alpha)},u)$ }

In this section, we determine the $\beta$-dual of the  spaces
$c_0(\Gamma,\Delta^{(\alpha)},u)$ and
$c(\Gamma,\Delta^{(\alpha)},u)$. For the sequence spaces $X$ and
$Y$, define the set $S(X,Y)$ by
\begin{equation}\label{3.1}
S(X,Y)=\{z=(z_{k})\in \omega: xz=(x_{k}z_{k})\in Y \textrm{ for all
} x\in X\}.
\end{equation}
With the notation of (\ref{3.1}), $\beta-$ dual of a sequence space
$X$ is defined by
$$
X^{\beta}=S(X,cs).
$$

%-------------------------------------------------------------------

\begin{lm}\label{e4.10}
$A\in (c_0:c)$ if and only if
\begin{equation}\label{e5.2}
\lim_{n\rightarrow\infty} a_{nk}=\alpha_k \ \textrm{for each fixed}\
k\in \mathbb{N},
\end{equation}

\begin{equation}\label{e5.3}
\sup_{n\in \mathbb{N}}\sum_k|a_{nk}|<\infty.
\end{equation}
\end{lm}

%-------------------------------------------------------------------

\begin{lm}\label{e4..2}
$A\in (c:c)$ if and only if (\ref{e5.2}) and (\ref{e5.3}) hold, and
\begin{equation}
\lim_{n\rightarrow\infty}\sum_ka_{nk}\quad\ \textrm{exists}.
\end{equation}
\end{lm}

\begin{lm}\label{lem3}
$A=(a_{nk})\in (\ell_\infty:\ell_\infty)$ if and only if
\begin{equation}\label{7c}
\sup_n\sum_{k} |a_{nk}|<\infty.
\end{equation}
\end{lm}

%-------------------------------------------------------------------

%-------------------------------------------------------------------

\begin{thm}\label{t3}
Define the sets $\Gamma_{1}$, $\Gamma_{2}$  and a matrix
$T=(t_{nk})$ by
$$
t_{nk}=\left\{\begin{array}{ll}
  \displaystyle t_{k}-t_{k+1},
   & (k<n)\\  \displaystyle t_{n}, & (k=n) \\\displaystyle\ 0, & (k>n)
\end{array}\right.
$$
for all $k,n\in\mathbb{N}$ where
$t_{k}=a_{k}\sum_{i=0}^{k}(-1)^i\frac{\Gamma(1-\alpha)}{i!\Gamma(1-\alpha-i)}
\frac{1}{u_{k-i}}$
\begin{eqnarray*}
\Gamma_{1}&=&\bigg\{ a=(a_n)\in
w:\lim_{n\rightarrow\infty}t_{nk}=\alpha_k ~\textrm{exists for
each}\ k\in \mathbb{N}\bigg\} \\
\Gamma_{2}&=&\bigg\{ a=(a_n)\in
w:\sup_{n\in\mathbb{N}}\sum_{k}^{}|t_{nk}|<\infty\bigg\}\\
\Gamma_{3}&=&\bigg\{ a=(a_n)\in
w:\sup_{n\in\mathbb{N}}\lim_{n\rightarrow\infty}\sum_k t_{nk}\quad\
\textrm{exists}\bigg\}.
\end{eqnarray*}

Then, $\{c_0(\Gamma,\Delta^{(\alpha)},u)\}^{\beta}=\Gamma_1\cap
\Gamma_2$ and
$\{c(\Gamma,\Delta^{(\alpha)},u)\}^{\beta}=\Gamma_1\cap \Gamma_2\cap
\Gamma_3$.

\end{thm}
\begin{pf}
We prove the theorem for the space
$c_0(\Gamma,\Delta^{(\alpha)},u)$. Let $a=(a_n)\in w$ and
$x=(x_k)\in c_0(\Gamma,\Delta^{(\alpha)},u) $. Then, we obtain the
equality
\begin{eqnarray}\label{beta}
\nonumber
\sum_{k=0}^{n}a_kx_k&=&\sum_{k=0}^{n}\bigg[a_{k}\sum_{i=0}^{\infty}(-1)^i\frac{\Gamma(1-\alpha)}{i!\Gamma(1-\alpha-i)}
\frac{1}{u_{k-i}}\bigg](y_{k-i}-y_{k-i-1})\\ \nonumber
&=&\sum_{k=0}^{n}\bigg[a_{k}\sum_{i=0}^{k}(-1)^i\frac{\Gamma(1-\alpha)}{i!\Gamma(1-\alpha-i)}
\frac{1}{u_{k-i}}\bigg](y_{k-i}-y_{k-i-1})\\ \nonumber
&=&\sum_{k=0}^{n-1}\bigg[a_{k}\sum_{i=0}^{k}(-1)^i\frac{\Gamma(1-\alpha)}{i!\Gamma(1-\alpha-i)}
\frac{1}{u_{k-i}}-a_{k+1}\sum_{i=0}^{k+1}(-1)^i\frac{\Gamma(1-\alpha)}{i!\Gamma(1-\alpha-i)}
\frac{1}{u_{k+1-i}}\bigg]y_k\\ \nonumber &
&+~\bigg[a_{n}\sum_{i=0}^{n}(-1)^i\frac{\Gamma(1-\alpha)}{i!\Gamma(1-\alpha-i)}
\frac{1}{u_{n-i}}\bigg]y_n\\  &=&T_ny
\end{eqnarray}

Then we deduce by (\ref{beta}) that $ax=(a_kx_k)\in cs$ whenever
$x=(x_k)\in c_0(\Gamma,\Delta^{(\alpha)},u)$ if and only if $T y\in
c$ whenever $y=(y_k)\in c_0$. This means that $a=(a_k)\in
\{c_0(\Gamma,\Delta^{(\alpha)},u)\}^\beta$ if and only if $T\in
(c_0:c)$. Therefore, by using Lemma \ref{e4.10}, we obtain :
\begin{eqnarray}
&&\lim_{n\rightarrow\infty}t_{nk}=\alpha_k\qquad \textrm{exists for
each}\
k\in\mathbb{N},\\
&&\sup_{n\in\mathbb{N}}\sum_{k}^{}|t_{nk}|<\infty.
\end{eqnarray}
Hence, we conclude that
 $\{c_0(\Gamma,\Delta^{(\alpha)},u)\}^{\beta}=\Gamma_1\cap
\Gamma_2$.

\end{pf}

%-------------------------------------------------------------------
\section{Some matrix transformations related to the sequence  spaces $c_0(\Gamma,\Delta^{(\alpha)},u)$
and $c(\Gamma,\Delta^{(\alpha)},u)$}

%-------------------------------------------------------------------
In this final section, we state some results which characterize
various matrix mappings on the spaces
$c_0(\Gamma,\Delta^{(\alpha)},u)$ and
$c(\Gamma,\Delta^{(\alpha)},u)$.
%-------------------------------------------------------------------

We shall write throughout for brevity that
\begin{equation}\label{z1}
\widetilde{a}_{nk}=z_{nk}-z_{n,k+1}~~\textrm{and}~~b_{nk}=\sum_{j=0}^{n}\bigg(\sum_{i=0}^{n-j}(-1)^i
\frac{\Gamma(\alpha+1)} {i!\Gamma(\alpha+1-i)}u_{i+j}\bigg) a_{jk}
\end{equation}
for all $k,n\in \mathbb{N}$  , where
$$z_{nk}=a_{nk}\sum_{i=0}^{k}(-1)^i\frac{\Gamma(1-\alpha)}{i!\Gamma(1-\alpha-i)}
\frac{1}{u_{k-i}}.$$

Now, we may give the following theorem.

%-------------------------------------------------------------------
\begin{thm}\label{teo4}
Let $\lambda$ be any given sequence space and $\mu\in \{c_{0},c\}$.
Then, $A =(a_{nk})\in (\mu(\Gamma,\Delta^{(\alpha)},u):\lambda) $ if
and only if ~$C\in (\mu: \lambda)$ and
\begin{eqnarray}\label{Bn}
C^{(n)}\in (\mu:c)
\end{eqnarray}
 for every fixed $n\in\mathbb{N}$,
where $c_{nk}=\widetilde{a}_{nk}$ and $C^{(n)}=(c_{mk}^{(n)})$ with
$$
c^{(n)}_{mk}=\left\{\begin{array}{ll}
  \displaystyle z_{nk}-z_{n,k+1},
   & (k<m)\\  \displaystyle z_{nm}, & (k=m) \\\displaystyle\ 0, & (k>m)
\end{array}\right.
$$
for all $k,m\in\mathbb{N}.$
\end{thm}
%-------------------------------------------------------------------
\begin{pf}
Let $\lambda$ be any given sequence space. Suppose that (\ref{z1})
holds between the entries of $A=(a_{nk})$ and $C=(c_{nk})$, and take
into account that the spaces $\mu(\Gamma,\Delta^{(\alpha)},u)$ and
$\mu$ are linearly isomorphic.

Let  $A =(a_{nk})\in (\mu(\Gamma,\Delta^{(\alpha)},u):\lambda) $ and
take any $y=(y_k)\in \mu$. Then, $C\Delta_u^{(\alpha)}(\Gamma)$
exists and $\{a_{nk}\}_{k\in \mathbb{N}}\in
\mu(\Gamma,\Delta^{(\alpha)},u)^{\beta}$ which yields that
(\ref{Bn}) is necessary and $\{c_{nk}\}_{k\in \mathbb{N}}\in
\mu^{\beta}$ for each $n\in \mathbb{N}$. Hence, $Cy$ exists for each
$y\in\mu$ and thus by letting $m\rightarrow\infty$ in the equality
\begin{equation}
\sum_{k=0}^{m}a_{nk}x_k=\sum^{m-1}_{k=0}
(z_{nk}-z_{n,k+1})y_{k}+z_{nm}y_m; \quad (m,n\in \mathbb{N})
\end{equation}
we have  that $Cy=Ax$ and so we have that $C\in (\mu:\lambda)$.

Conversely, suppose that $C\in (\mu:\lambda)$ and (\ref{Bn}) hold,
and take any $x=(x_{k})\in \mu(\Gamma,\Delta^{(\alpha)},u)$. Then,
we have $\{c_{nk}\}_{k\in \mathbb{N}}\in \mu^{\beta}$ which gives
together with (\ref{Bn}) that $\{a_{nk}\}_{k\in \mathbb{N}}\in
\mu(\Gamma,\Delta^{(\alpha)},u)^{\beta}$ for each $n\in \mathbb{N}$.
So, $Ax$ exists. Therefore, we obtain from the equality
\begin{equation}
\sum_{k=0}^{m}
c_{nk}y_{k}=\sum_{k=0}^{m}\bigg[\sum_{j=k}^{m}\bigg(\sum_{i=0}^{j-k}(-1)^i\frac{\Gamma(\alpha+1)}
{i!\Gamma(\alpha+1-i)}u_{i+k}\bigg) c_{nj}\bigg]x_{k} \quad
\textrm{for all}\  n\in \mathbb{N},
\end{equation}
as $m\rightarrow \infty$ that $Ax=Cy$ and this shows that $A\in
(\mu(\Gamma,\Delta^{(\alpha)},u):\lambda)$. This completes the
proof.
\end{pf}

%%-------------------------------------------------------------------
%By changing the roles of the spaces
%$\mu(\Gamma,\Delta^{(\alpha)},u)$ and $\mu$ with $\lambda$, we have
%%-------------------------------------------------------------------
\begin{thm}\label{teo5}
Suppose that the entries of the infinite matrices $A=(a_{nk})$ and
$B=(b_{nk})$ are connected with the relation (\ref{z1}) and
$\lambda$ be given sequence space and $\mu\in \{c_0,c\}$. Then
$A=(a_{nk})\in(\lambda:\mu(\Gamma,\Delta^{(\alpha)},u))$ if and only
if $B=(b_{nk})\in(\lambda:\mu)$ .
\end{thm}

%-------------------------------------------------------------------
\begin{pf}
Let $z=(z_k)\in\lambda$ and consider the following equality with
(\ref{z1})
$$
\sum_{k=0}^{m}
b_{nk}z_{k}=\sum_{j=0}^{n}\bigg[\sum_{k=0}^{m}\bigg(\sum_{i=0}^{n-j}(-1)^i\frac{\Gamma(\alpha+1)}
{i!\Gamma(\alpha+1-i)}u_{i+j}\bigg) a_{jk}\bigg]z_{k} \quad \
(m,n\in \mathbb{N}),
$$
which yields as $m\rightarrow\infty$ that
$(Bz)_n=[\Delta_u^{(\alpha)}(\Gamma)(Az)]_n$. Hence, we obtain that
$Az\in\mu(\Gamma,\Delta^{(\alpha)},u)$ whenever $z\in\lambda$ if and
only if $Bz\in\mu$ whenever $z\in\lambda$.
\end{pf}

%-------------------------------------------------------------------
 We will have several consequences by using Theorem \ref{teo4} and Theorem
 \ref{teo5}. But we must give firstly some
 relations which are important for consequences:

%-------------------------------------------------------------------
\begin{eqnarray}
\label{co1}&& \sup_{n\in\mathbb{N}}\sum_k|a_{nk}|<\infty \\
\label{co2}&&  \lim_{n\rightarrow\infty}a_{nk}=\alpha_k\quad
\textrm{exists for each fixed}\ k\in \mathbb{N}
\\
\label{co3}&& \lim_{k\rightarrow\infty}a_{nk}=0 \quad \textrm{for
each fixed}\ n\in \mathbb{N}\\
\label{co4}&& \lim_{n\rightarrow\infty}\sum_ka_{nk}\quad\ \textrm{exists}\\
\label{co5}&&\lim_{n\rightarrow\infty}\sum_ka_{nk}=0\\
\label{co6}&& \sup_{K\in \mathcal{F}}\sum_n\bigg|\sum_{k\in
K}a_{nk}\bigg|<\infty\\
\label{co7}&& \lim_{n\rightarrow\infty}\sum_k |a_{nk}|=0\\
\label{co8}&&\sup_{n,k}|a_{nk}|<\infty\\
\label{co9}&&\lim_{m\rightarrow\infty}\sum_k|a_{nk}|=\sum_k|\alpha_{k}|
\end{eqnarray}

%-------------------------------------------------------------------

Now, we can give the corollaries:
%------------------------------------------------------------------------------
\begin{cor}
The following statements hold:

(i) $A=(a_{nk})\in (c_0(\Gamma,\Delta^{(\alpha)},u):\ell_\infty)=
(c(\Gamma,\Delta^{(\alpha)},u):\ell_\infty)$ if and only if
(\ref{co1}) holds with $\tilde{a}_{nk}$ instead of $a_{nk}$ and
(\ref{Bn}) also holds.

(ii) $A=(a_{nk})\in (c_0(\Gamma,\Delta^{(\alpha)},u):c)$ if and only
if (\ref{co1})and (\ref{co2}) hold with $\tilde{a}_{nk}$ instead of
$a_{nk}$ and (\ref{Bn}) also holds.

(iii) $A=(a_{nk})\in (c_0(\Gamma,\Delta^{(\alpha)},u):c_0)$ if and
only if (\ref{co1}) and (\ref{co3})  hold with $\tilde{a}_{nk}$
instead of $a_{nk}$ and (\ref{Bn}) also holds.

(iv) $A=(a_{nk})\in (c(\Gamma,\Delta^{(\alpha)},u):c)$ if and only
if (\ref{co1}), (\ref{co2}) and (\ref{co4}) hold with
$\tilde{a}_{nk}$ instead of $a_{nk}$ and (\ref{Bn}) also holds.

(v) $A=(a_{nk})\in (c(\Gamma,\Delta^{(\alpha)},u):c_0)$ if and only
if (\ref{co1}), (\ref{co3}) and (\ref{co5}) hold with
$\tilde{a}_{nk}$ instead of $a_{nk}$ and (\ref{Bn}) also holds.

(vi) $A=(a_{nk})\in
(c_0(\Gamma,\Delta^{(\alpha)},u):\ell)=(c(\Gamma,\Delta^{(\alpha)},u):\ell)$
if and only if (\ref{co6}) holds with $\tilde{a}_{nk}$ instead of
$a_{nk}$ and (\ref{Bn}) also holds.
\end{cor}
%---------------------------------------------------------------------------------------
\begin{cor}
The following statements hold:

(i)$A=(a_{nk})\in (\ell_\infty:c_0(\Gamma,\Delta^{(\alpha)},u))$ if
and only if (\ref{co7}) hold with $b_{nk}$ instead of $a_{nk}$.

(ii)$A=(a_{nk})\in (c:c_0(\Gamma,\Delta^{(\alpha)},u))$ if and only
if (\ref{co1}), (\ref{co3}) and (\ref{co5}) hold with $b_{nk}$
instead of $a_{nk}$.

(iii) $A=(a_{nk})\in (c_0:c_0(\Gamma,\Delta^{(\alpha)},u))$ if and
only if (\ref{co1}) and (\ref{co3})  hold with $b_{nk}$ instead of
$a_{nk}$.

(iv)$A=(a_{nk})\in (c:c_0(\Gamma,\Delta^{(\alpha)},u))$ if and only
if  (\ref{co3}) and (\ref{co8}) hold with $b_{nk}$ instead of
$a_{nk}$.

(v) $A=(a_{nk})\in (\ell_\infty:c(\Gamma,\Delta^{(\alpha)},u))$ if
and only if (\ref{co2}) and (\ref{co9}) hold with $b_{nk}$ instead
of $a_{nk}$.

(vi) $A=(a_{nk})\in (c:c(\Gamma,\Delta^{(\alpha)},u))$ if and only
if (\ref{co1}), (\ref{co2}) and (\ref{co4}) hold with $b_{nk}$
instead of $a_{nk}$.

(vii) $A=(a_{nk})\in (c_0:c(\Gamma,\Delta^{(\alpha)},u))$ if and
only if (\ref{co1}) and (\ref{co2}) hold with $b_{nk}$ instead of
$a_{nk}$.

(viii) $A=(a_{nk})\in (\ell:c(\Gamma,\Delta^{(\alpha)},u))$ if and
only if (\ref{co2}) and (\ref{co8}) hold with $b_{nk}$ instead of
$a_{nk}$.

\end{cor}

\begin{thebibliography}{30}
%-----------------------------------------------------------------

\bibitem{ggl} G.G. Lorentz, A contribution to the theory of
divergent sequences, Acta Math. \textbf{80}(1948), 167-190.

\bibitem{kzmaz} H. K{\i}zmaz,On certain sequence spaces, Canad. Math.
bull., 24(1981) 169-176

\bibitem{met} M. Et, On some difference sequence spaces, Turkish J.
Math. 17 (1993) 18-24.

\bibitem{metolak} M. Et, R. \c{C}olak, On some generalized difference sequence
spaces, Soochow J. Math. 21 (1995) 377-386.

\bibitem{cem} R. \c{C}olak, M. Et, E. Malkowsky, Some topics of sequence
spaces, in: Lecture Notes in Mathematics, Fýrat Univ. Press, Turkey,
2004, pp.1-63. Fýrat Univ. Elazýð, 2004, ISBN:975-394-038-6.

\bibitem{w} C. S. Wang, On N\"{o}rlund sequence spaces, Tamkang J.
Math. \textbf{9}(1978), 269-274.

\bibitem{nglee} P. N. Ng, P.Y. Lee, Ces\`{a}ro sequence spaces of
non-absolute type, Comment. Math. Prace Mat. \textbf{20}(2)(1978),
429-433.

\bibitem{ajem} A. M. Jarrah, E. Malkowsky, BK-spaces, bases and
linear operators, Ren. Circ. Mat. Palermo II \textbf{52}(1990),
177-191.

\bibitem{emmmss} E. Malkowsky, M. Mursaleen, S. Suantai, The dual
spaces of sets of difference sequences of order $m$ and matrix
transformations, Acta Math. Sin. Eng. Ser. 23(3)(2007) 521-532.

\bibitem{ba} B. Altay, On the space of $p-$ summable difference
sequences of order $m$, $(1\leq p<\infty)$, Studia Sci. Math.
Hungar. 43(4)(2006) 387-402.

\bibitem{batfbr} B. Altay, F. Ba\c{s}ar, Generalization of the space $\ell(p)$ derived by weighted mean ,
J. Math. Anal. Appl 330 (2007) 174-185.

\bibitem{bafb5} B. Altay, F. Ba\c{s}ar, Certain topological properties and duals
of the matrix domain of a triangle matrix in a sequence space, J.
Math. Anal. Appl. 336(1)(2007) 632-645.

\bibitem{mbeek1} M. Ba\c{s}ar{\i}r, E. E. Kara, On the $m^{th}$
order difference sequence space of generalized weighted mean and
compact operators, Acta Mathematica Scientia, 33(B3)(2013), 1-18.

\bibitem{mbeek3} M. Ba\c{s}ar{\i}r, E. E. Kara, On some
difference sequence space of  weighted means and compact operators,
Annals of Functional Analysis, 2(2)(2011), 116-131.

\bibitem{sdcc} S. Demiriz, C. \c{C}akan, Some topological and geometrical properties of a new
difference sequence space, Abstract and Applied Analysis, In Press.

\bibitem{od} Duyar O., Demiriz, S., On some new generalized difference sequence spaces and
their topological properties, Journal of New Result in Science,
6(2014) 1-14.


\bibitem{cafb} C. Ayd{\i}n, F.Ba\c{s}ar, Some new difference sequence spaces,
Applied Mathematics and Computation, \textbf{157}(2004), 677-693.

\bibitem {bail1} P. Baliarsingh, Some new difference sequence spaces of fractional
order and their dual spaces, Applied Mathematics and Computation,
219(2013) 9737-9742.

\bibitem {bail2} S. Dutta, P. Baliarsingh,  A note on paranormed
difference sequence spaces of fractional order and their matrix
transformations, Journal of the Egyptian Mathematical Society,
(2013).

\bibitem {bail3} P. Baliarsingh, S. Dutta, A unifying approach to
the difference operators and their applications, Bol. Soc. Paran.
Mat., (3s)\textbf{v. 33} 1 (2015): 49-57.

\bibitem{fb} F. Ba\c{s}ar, Summability Theory and Its Appliactions,
Bentham Science Publishers, ISBN:978-1-60805-252-3, 2011.

\bibitem{msht} M. Stieglitz, H. Tietz, Matrix transformationen von
folgenr\"{a}umen eine ergebnisübersicht, Math. Z. 154(1977) 1-16.

\bibitem{m} I.J. Maddox, Elements of Functional Analysis, second ed., The Cambridge
University Press, Cambridge, 1988.

\bibitem{aw} A. Wilansky, Summability through Functional Analysis,
North-Holland Mathematics Studies , Amsterdam, \textbf{85} 1984.




%-----------------------------------------------------------------
\end{thebibliography}
\end{document}